\documentclass[reqno,11pt]{amsart}
\usepackage{amsmath}
\usepackage{graphicx, color, enumerate}
\usepackage{amscd}
\usepackage{amsmath}
\usepackage{amsfonts}
\usepackage{amssymb}
\usepackage{mathrsfs}
\usepackage{latexsym}

\newcommand{\be}{\begin{equation}}
\newcommand{\ee}{\end{equation}}
\newenvironment{pf}{\noindent{\it
Proof}.\enspace}{\rule{2mm}{2mm}\medskip}

\newcommand{\R}{\mathbb{R}}
 \newcommand{\Rn}{\mathbb{R}^n}

\newcommand{\E}{\mathbb{E}}

\newcommand{\h}{\mathbb{H}}
\newcommand{\dyle}{\displaystyle}
\renewcommand{\a }{\alpha }
\renewcommand{\b }{\beta }

\newcommand{\D }{\Delta }
\newcommand{\e }{\varepsilon }
\newcommand{\g }{\gamma}
\newcommand{\G }{\Gamma }
\renewcommand{\l }{\lambda }
\renewcommand{\L }{\Lambda }
\newcommand{\m }{\mu }
\newcommand{\n }{\nabla }

\renewcommand{\o }{\omega }

\newcommand{\cN}{{\mathcal{N}}}

\newcommand{\intn}{\int_{\Rn}}
\newcommand{\intR}{\int_\R}

\newcommand{\wt}{\widetilde}

\newcommand{\bu}{{\bf u}}
\newcommand{\bv}{{\bf v}}

\newcommand{\bw}{{\bf w}}
\newcommand{\bo}{{\bf 0}}

\newcommand{\bh}{{\bf h}}

\newtheorem{Theorem}{Theorem}

\newtheorem{Lemma}[Theorem]{Lemma}
\newtheorem{Proposition}[Theorem]{Proposition}
\newtheorem{Definition}[Theorem]{Definition}
\newtheorem{remark}[Theorem]{Remark}
\newtheorem{remarks}[Theorem]{Remarks}
\newtheorem{example}[Theorem]{Example}
\newtheorem{examples}[Theorem]{Examples}
\newenvironment{Remark}{\begin{remark} \rm}{\rule{2mm}{2mm}\end{remark}}
\newenvironment{Remarks}{\begin{remarks}
\rm}{\rule{2mm}{2mm}\end{remarks}}

        \usepackage{latexsym}
        \usepackage{amsmath}
        \usepackage{amsfonts}
        \usepackage{amssymb}
\author[Eduardo Colorado]{Eduardo Colorado$^{*}$}
\address{\noindent Departamento de Matem\'aticas, Universidad Carlos III de Madrid,
Avenida de la Universidad 30, 28911 Legan\'es, Madrid, Spain $\&$
Instituto de Ciencias Matem\'aticas, ICMAT (CSIC-UAM-UC3M-UCM),
 C/Nicol\'as Cabrera 15, 28049 Madrid, Spain.}
\email{eduardo.colorado@uc3m.es, eduardo.colorado@icmat.es}
\thanks{$^*$Partially supported by Ministry of Economy and Competitiveness of
Spain and FEDER funds under research project MTM2013-44123-P}

\title[Existence of bound and ground states for some coupled NLS-KdV equations]{On the existence of
bound and ground states for some coupled nonlinear Schr\"odinger--Korteweg-de Vries equations}

\begin{document}
{\sf \maketitle

\noindent {\small  {\bf Abstract.} We demonstrate existence of
positive bound and ground states for a system of coupled nonlinear
Schr\"odinger--Korteweg-de Vries  equations. More precisely, we
prove there is a positive radially symmetric ground state if either
the coupling coefficient $\b>\L$ (for an appropriate constant $\L>0$) or
$\b>0$ with appropriate conditions on the other parameters of the
problem. Concerning bound states, we prove there exists a positive
radially symmetric bound state if either $0<\b$ is sufficiently
small or $0<\b<\L$ with some appropriate conditions on the
parameters as for the ground states. That results give a
classification of positive solutions as well as multiplicity of
positive solutions. Furthermore, we consider a system with more
general power nonlinearities, proving the above results, and also we
study natural extended systems with more than two equations.
Although the techniques we employed are variational, we look for
critical points of an energy functional different from the classical
one used in this kind of systems. Our approach improves many of the
previous known results, as well as permit us to show new results not
considered or studied before.}

\

\tableofcontents
%%%%%%%%%%%%%%%%%%%%%%%%%%%%%%
\section{Introduction}\label{sec:intro}
%%%%%%%%%%%%%%%%%%%%%%%%%%%%%%
The aim of this work is to study a system of coupled  nonlinear
Schr\"odinger--Korteweg-de Vries (NLS-KdV for short) equations as
follows,
\begin{equation}\label{NLS-KdV}
\left\{\begin{array}{rcl}
if_t + f_{xx} + |f|^2f& = & \b fg\\
g_t +g_{xxx} +gg_x   & = & \frac 12\b(|f|^2)_x
\end{array}\right.
\end{equation}
where $f=f(x,t)\in \mathbb{C}$ while $g=g(x,t)\in \mathbb{R}$, and
$\b>0$ is the real coupling coefficient. System \eqref{NLS-KdV}
appears in phenomena of interactions between short and long
dispersive waves, arising in fluid mechanics, such as  the
interactions of capillary - gravity water waves. Indeed, $f$
represents the short-wave, while $g$ stands for the long-wave. See
\cite{aa,cl,fo} and the references therein for more details.

We  look for  solitary ``traveling-waves", namely solutions to
\eqref{NLS-KdV} of the form \be\label{eq:sol-t-v-solutions}
(f(x,t),g(x,t))=\left(e^{i\o t} e^{i\frac c2
x}u(x-ct),v(x-ct)\right)\quad\mbox{with}\quad u,\, v\quad\mbox{real
functions.} \ee Choosing  $\l_1=\o+\frac{c^2}{4}$, $\l_2=c$, we get
that $u,\, v$ solve the following problem
\begin{equation}\label{NLS-KdV2}
\left\{\begin{array}{rcl}
-u'' +\l_1 u & = & u^3+\beta uv \\
-v'' +\l_2 v & = & \frac 12 v^2+\frac 12\beta u^2.
\end{array}\right.
\end{equation}
This system has been previously studied by Dias et al. in
\cite{dfo}. Also a generalization of \eqref{NLS-KdV2}, see
\eqref{NLS-KdV-pq-general2}, has been previously analyzed  by the
same authors in \cite{dfo2} and by Albert and Bhattarai in
\cite{ab}. A comparison with our results will be done at the end of
this introduction.

\

The main goal of this work  is three fold. One is to give a
classification of positive solutions of  \eqref{NLS-KdV2}, proving:

\noindent -Existence of positive even ground states of \eqref{NLS-KdV2} under the
following hypotheses:
\begin{itemize}
\item  the coupling coefficient $\b>\L>0$ for an appropriate constant $\L$; see Theorem \ref{th:1},
\item $\b>0$ and $\l_2\gg 1$; see Theorem \ref{th:ground2}.
\end{itemize}
-Existence of positive even bound states of \eqref{NLS-KdV2} when:
\begin{itemize}
\item $0<\b\ll 1$; see Theorem \ref{th:2}, where we also give a bifurcation result,
\item $0<\b<\L$ and $\l_2\gg 1$; see Theorem \ref{th:ground2}.
\end{itemize}
The coexistence of positive bound and ground states for  $0<\b<\L$
and $\l_2$ large is a great novelty and difference with the more
studied systems of NLS equations in the last several years; see
Remark \ref{rem:11}-$(ii)$.

The second goal is that we study a more general system than
\eqref{NLS-KdV2}, with more general power nonlinearities given by
\eqref{NLS-KdV-pq-general2}, for which we show that previous results
for \eqref{NLS-KdV2} hold with similar conditions on the
coefficients.

We also analyzed a particular case of $\l_1,\,\l_2$ in which there
exists an explicit positive solution.

The last goal is to consider natural extensions of \eqref{NLS-KdV2}
to systems with more than two equations, as well as deal with
extensions to the dimensional cases $n=2,3$, for which although
\eqref{NLS-KdV} has no sense, the stationary system \eqref{NLS-KdV2}
makes sense and can be seen, for example, as the stationary system
when one looks for standing wave solutions of the corresponding
evolutionary system of NLS equations. For some of these extended
problems we show similar results as described above on the existence
of positive radially symmetric  bound and ground state solutions.
Other systems with at least two NLS equations and at least one KdV
equation will be analyzed (in more detail) in a forthcoming paper.

\

Besides the previous achievements, it is relevant to point out that is the first
time that  our variational procedure (in part developed in \cite{ac1,ac2} for systems of coupled NLS equations)
is employed to study coupled NLS-KdV equations in an appropriate way, see Remark \ref{rem:compacidad}-$(ii)$. Even
more, it seems to be better in many ways than the classical approach
used before to study NLS-KdV systems as we will see. Note that our
method could be exploited to study related problems.

\

 It is worth pointing out that,  for any $\b\in \R$,  System
\eqref{NLS-KdV2}  has a unique {\it semi-trivial} positive solution
$\bv_2=(0,V_2)$, where
$V_2(x)=3\,\l_2\,\mbox{sech}^2\!\left(\frac{\sqrt{\l_2}}{2}x\right)$
is the unique positive even solution of  $-v''+\l_2 v=\frac 12 v^2$
in $W^{1,2}(\R)$; \cite{kw}. As might be expected, we look for
different solutions from the preceding one. We are interested not in
non-negative solutions but  in positive ones, and therefore
different from $\bv_2$.

As we announced above, a comparison with our results and the
previous works \cite{dfo,dfo2,ab} is in order. In \cite{dfo}, Dias
et all. studied \eqref{NLS-KdV2} in the particular case $\l_1=\l_2$
and they proved the existence of non-negative bound state solutions
when the coupling parameter $\b>\frac 12$. Here, we have improved
that result in three ways. First, we have considered $\l_1$ not
necessarily equals $\l_2$ and we proved not only the existence of
non-negative bound states but also positive even ground states for
$\b$ greater than a constant $\L>0$ defined by \eqref{eq:Lambda},
for which in the setting of \cite{dfo}, we have $\L\le \frac 12$.
Secondly, we show the existence of positive even bound states when
$0<\b\ll1$, not studied in \cite{dfo}. Thirdly, we also show that if
$\l_2$ is sufficiently large, there exists a positive even ground
state for every $\b>0$ (with $\l_1$ not necessarily equals $\l_2$)
and a positive even bound state provided $0<\b<\L$. In \cite{dfo2},
among other results, Dias et all. studied System
\eqref{NLS-KdV-pq-general2} with $2<q<5$, $p\in\{ 2,3,4\}$,
$\m_2=p+1$, and they established (\cite[Theorem 4.1]{dfo2}) the
existence of a non-trivial bound state solution for all $\b>0$ if
$p=3,4$ and $\b>3$ if $p=2$. Finally, in \cite{ab} Albert and
Bhattarai studied, among other topics, system
\eqref{NLS-KdV-pq-general2}  in a more general setting than in
\cite{dfo2}, precisely they considered $2\le q<5$, $2\le p< 5$ with
$p$ a rational number with odd denominator; they proved  the
existence of a positive even bound state for each $\b>0$, improving
the above cited  result by \cite{dfo2}. In our manuscript we consider
$2\le p<\infty$, $2\le q<\infty$ and we prove that there exists a
positive even ground state of \eqref{NLS-KdV-pq-general2} if either
$\b>\L$ or $\b>0$, $q>2p-2$ and $\l_2$ is large enough. Concerning
bound states we show the existence of a positive even bound state of
\eqref{NLS-KdV-pq-general2} if either $0<\b<\L$, $q>2p-2$ and $\l_2$
is large enough
 or $0<\b\ll 1$,  proving also in this last case a
bifurcation result for the bound state we find.

By the discussion above of problem \eqref{NLS-KdV-pq-general2},
 we improve and extend some of the results by \cite{dfo2,ab}.
Additionally, we establish some new results for
\eqref{NLS-KdV-pq-general2}, as the multiplicity one of coexistence
of positive bound and ground states in some range of the parameters.
Finally, we analyzed extended systems of \eqref{NLS-KdV2} with more
than two equations, that up to our knowledge, have not been
considered previously in the literature. Also, we study the
qualitative and quantitative properties of the explicit solutions of
\eqref{NLS-KdV2}.

\

A preliminary announcement of some results in the present work appeared in \cite{c2}.

\

The paper is organized as follows. In Section \ref{sec:funct} we
introduce the functional framework, notation and give some
definition. Next,  we define the Nehari Manifold in Section
\ref{sec:key}, proving some properties of it, we establish a useful
measure lemma and show a result dealing with qualitative properties
of the semi-trivial solution. Section \ref{sec:ground} is divided
into two subsections, the first one contains the existence of ground
states, and the second one deals with the existence of bound states.
In Section \ref{sec:extended} we study a  system with more general
power nonlinearities, proving similar results as in the previous one
with the appropriate changes. Section \ref{sec:remarks}
contains two subsections, the first one deals with an explicit
solution, while the last one is devoted to study natural extensions
to systems with more than two equations.
%%%%%%%%%%%%%%%%%%%%%%%%%%%%%%%%%%%%
\section{Functional setting and notation}\label{sec:funct}
%%%%%%%%%%%%%%%%%%%%%%%%%%%%%%%%%%%%
Let $E$ denotes the Sobolev space $W^{1,2}(\R)$, that can be defined
as the completion of $\mathcal{C}_0^1(\R)$ endowed with the norm
$$
\|u\| =\sqrt{(u\mid u)},
$$
which comes from the scalar product
$$
(u\mid w)=\intR (u'w'+uw)dx.
$$
We will denote the following equivalent norms and scalar products in $E$,
$$
\|u\|_j=\| u\|_{\l_j}=\left(\intR (| u'|^2+\l_j u^2)\,
dx\right)^{\frac 12},
$$
$$
( u|v)_j=(u|v)_{\l_j}=\intR ( u'\cdot v'+ \l_j uv)\, dx; \quad
j=1,2.
$$
Let us define the product Sobolev space $\E=E\times E$. The elements in $\E$ will be denoted by $\bu =(u,v)$, and $\bo=(0,0)$.
We will take
$$
\|\bu\|=\sqrt{\|u\|_1^2+\|v\|_2^2}
$$
as a norm in $\E$.

For $\bu=(u,v)\in \E$, the notation $\bu\geq \bo$, resp. $\bu>\bo$,
means that $u,v\geq 0$, resp. $u,v>0$. We denote $H$ as the space of
even (radially symmetric) functions in $E$, and $\h=H\times H$.

We define the functionals
$$
I_1(u)=\tfrac 12 \|u\|_1^2 -\tfrac 14\, \intR u^4dx,\qquad I_2(v)=\tfrac 12 \|v\|_2^2 -\tfrac 16\, \intR v^3dx,\qquad u,\, v\in E,
$$
and
$$
\Phi (\bu)= I_1(u)+I_2(v)- \tfrac 12\b \intR u^2v\,dx,\qquad \bu\in \E.
$$
We also write
$$
G_\b (\bu)=\tfrac 14\, \intR u^4dx+\tfrac 16\, \intR v^3dx+ \tfrac 12\b \intR u^2v\,dx, \qquad \bu\in \E,
$$
and using this notation we can rewrite
the energy functional
$$
\Phi (\bu)=\frac 12\|\bu \|^2 -G_\b (\bu),\quad \bu\in \E.
$$
\begin{Definition}
We say that $\bu\in \E$ is a non-trivial {\it bound state} of \eqref{NLS-KdV2} if $\bu$ is a non-trivial  critical point
of $\Phi$.
A bound state $\wt{\bu}$ is called ground state if its energy is minimal among all the non-trivial bound states, namely
\begin{equation}\label{eq:gr}
\Phi(\wt{\bu})=\min\{\Phi(\bu): \bu\in \E\setminus\{\bo\},\; \Phi'(\bu)=0\}.
\end{equation}
\end{Definition}
%%%%%%%%%%%%%%%%%%%%%%%%%%%%%%
\section{Nehari manifold and key results}\label{sec:key}
%%%%%%%%%%%%%%%%%%%%%%%%%%%%%%
We will work mainly  in  $\h$.
Setting
$$
\Psi(\bu)=(\n\Phi(\bu)|\bu)=(I_1'(u)|u)+(I_2'(v)|v)-\frac 32\b \intR u^2v\,dx,$$
we define the corresponding Nehari manifold
$$
\cN =\{ \bu\in \h\setminus\{\bo\}: \Psi (\bu)=0\}.
$$
Then, one has that
\begin{equation}\label{eq:gamma}
(\n \Psi(\bu) \mid \bu)= - \|\bu \|^2-\intR u^4\,dx<0,\quad\forall\, \bu\in \cN,
\end{equation}
thus $\cN$ is a smooth manifold locally near any point $\bu\not= \bo$ with $\Psi(\bu)=0$. Moreover,
$\Phi''(\bo)= I_1''(0)+I_2''(0)$ is positive definite, so we infer that $\bo$ is a strict minimum for $\Phi$. As a consequence,  $\bo$ is an isolated point of the
set $\{\Psi(\bu)=0\}$, proving that, on the one hand $\cN$ is a smooth complete manifold of codimension $1$, and on the other hand there exists a constant $\rho>0$ so that
\be\label{eq:bound}
\|\bu\|^2>\rho,\qquad\forall\,\bu\in \cN.
\ee
Furthermore, \eqref{eq:gamma} and \eqref{eq:bound} plainly imply that $\bu\in \h\setminus\{\bo\}$ is a critical point of $\Phi$ if and only if $\bu\in\cN$ is a critical point of $\Phi$ constrained on $\cN$.

\begin{remarks}\label{rem:obs1}
\begin{itemize}
\item[(i)] By the previous arguments, the Nehari manifold $\cN$ is a natural constraint of $\Phi$. Also, it is relevant to point out that working on the Nehari manifold,
the functional $\Phi$ takes the form:
\be\label{eq:restriction0}
\Phi|_{\cN}(\bu)= \frac 16\|\bu\|^2+\frac{1}{12}\intn u^4dx=:F(\bu),
\ee
and by using  \eqref{eq:bound} into \eqref{eq:restriction0} we have
\begin{equation}\label{eq:restriction}
\Phi|(\bu)\ge   \frac 16\|\bu\|^2>\frac 16 \rho\qquad \forall\, \bu\in\cN.
\end{equation}
Therefore, \eqref{eq:restriction} shows that the functional $\Phi$ is bounded from below  on $\cN$, so one can try to minimize it on the Nehari manifold.
\item[(ii)] With respect to the Palais-Smale (PS for short) condition, we recall that in the one dimensional case,
one cannot expect a compact embedding of $E$ into $L^q(\R)$ for $2<
q<\infty$. Indeed, working on $H$ (the radial or even case) is not
true too; see \cite[Remarque I.1]{Lions-JFA82}. However, we will
show that for a PS sequence we can find a subsequence for which the
weak limit is a solution. This fact jointly with some properties of
the Schwarz symmetrization will permit us to prove the existence of
positive even ground states in Theorem \ref{th:1}.
\end{itemize}
\end{remarks}
Due to the lack of compactness mentioned above in Remarks \ref{rem:obs1}-$(ii)$, we state a measure theory result given in \cite{lions2}  that we will use in the proof of Theorem \ref{th:1}.
\begin{Lemma}\label{lem:measure}
If $2<q<\infty$, there exists a constant $C>0$ so that
\begin{equation}\label{eq:measure}
\intn |u|^q \, dx\le C\left( \sup_{z\in\R}\int_{|x-z|<1}|u(x)|^2dx\right)^{\frac{q-2}{2}}
\| u\|^2_{E},\quad \forall\: u\in E.
\end{equation}
\end{Lemma}
See \cite{c-fract} for an extension of this lemma in fractional Sobolev spaces, and an application of it in a fractional system of NLS equations.

\

Let $V$ denotes the unique positive even solution of $-v''+v=v^2$,
$v\in H$; see \cite{kw}. Setting \be\label{eq:segunda}
V_2(x)=2\l_2\,V(\sqrt{\l_2}\,x)=\frac{3\l_2}{\cosh^2\left(\frac{\sqrt{\l_2}}{2}x\right)},
\ee one has that $V_2$ is the unique positive  solution of
$-v''+\l_2v=\frac 12 v^2$ in $H$. Hence $\bv_2 := (0,V_2)$ is a
particular solution of  \eqref{NLS-KdV2} for any $\b\in\R$, and
moreover, it is the unique non-negative semi-trivial solution of
\eqref{NLS-KdV2}. We also define the corresponding Nehari manifold,
$$
\cN_2 =\left\{v\in H : (I_2'(v)|v)=0\right\}=\left\{v\in H : \|v\|_2^2 -\frac 12\intR v^3dx=0\right\}.
$$
Let us denote $T_{\bv_2}\cN$ the tangent space  to  $\cN$ on $\bv_2$. Since
$$
\bh=(h_1,h_2)\in  T_{\bv_2}\cN \Longleftrightarrow
(V_2|h_2)= \frac 34\intR V_2^2h_2\,dx ,
$$
it follows that
\begin{equation}\label{eq:tang1}
(h_1,h_2)\in T_{\bv_2} \cN  \Longleftrightarrow h_2\in T_{V_2} \cN_2.
\end{equation}

\begin{Proposition}\label{lem:gl3}
There exists $\L>0$ such that:
\begin{itemize}
\item[(i)]  if $\b< \L$, then $\bv_2$ is a strict local minimum  of $\Phi$ constrained on $\cN$,
\item[(ii)] for any   $\b>\L$, then $\bv_2$ is a  saddle point of $\Phi$ constrained  on $\cN$. Moreover,  $\dyle\inf_{\cN}\Phi <\Phi(\bv_2)$.
\end{itemize}
\end{Proposition}
\begin{pf} $(i)$ We define
\be\label{eq:Lambda} \L=\inf_{\varphi\in H\setminus\{
0\}}\frac{\|\varphi\|_1^2}{\intR V_2\varphi^2dx}. \ee One has that
for $\bh\in  T_{\bv_2}\cN$, \be\label{eq:Phi-segunda}
\Phi''(\bv_2)[\bh]^2 =\|h_1\|_1^2 +I_2''(V_2)[h_2]^2-\b\intR V_2
h_1^2dx. \ee Let us take $\bh=(h_1,h_2)\in T_{\bv_2}\cN$, by
\eqref{eq:tang1} $h_2\in T_{V_2} \cN_2$, then using that $V_2$ is
the minimum of $I_2$ on $\cN_2$, there exists a  constant $c>0$ so
that \be\label{eq:minimo-pos} I_2'' (V_2)[h_2]^2\ge c\|h_2\|_2^2.
\ee Let $h_1$ be a function with $\| h_1\|_1^2=\L \intR
V_2h_1^2\,dx$, i.e., which exists since the infimum $\L$ defined by
\eqref{eq:Lambda} is achieved, using this fact jointly with
\eqref{eq:Phi-segunda} and $\b<\L$, there exists another constant
$c_1>0$ so that, \be\label{minimo} \Phi''(\bv_2)[\bh]^2 \ge
c_1\|h_1\|_1^2 +c_2\| h_2\|^2. \ee Notice that $\Phi'(\bv_2)=0$
implies that $D^2\Phi_{\cN}(\bv_2)[\bh]^2=\Phi''(\bv_2)[\bh]^2$ for
all $\bh\in T_{\bv_2}\cN$, and thus using \eqref{minimo} we infer
that  $\bv_2$ is a local strict minimum of $\Phi$ on $\cN$.

\

$(ii)$
According to \eqref{eq:tang1}, $\bh=(h_1,0)\in T_{\bv_2}\cN$ for any $h_1\in H$.
We have that, for $\b>\L$, there exists $\wt{h}\in H$ with
$$
\L< \frac{\|\wt{h}\|_1^2}{\intR V_2\wt{h}^2dx}<\b,
$$
thus, taking $\bh_0=(\wt{h},0)\in T_{\bv_2}\cN$, by \eqref{eq:Phi-segunda} we find
$$
\Phi''(\bv_2)[\bh_0]^2 =\|\wt{h}\|_1^2 -\b\intR V_2 \wt{h}^2dx<0,
$$
finishing the proof taking $\L=\L$. \end{pf}
\begin{remark}\label{rem:1}
 If one consider $\l_1=\l_2$ as in \cite{dfo}, taking $\bh_0=(V_2,0)\in T_{\bv_2}\cN$ in the proof of Proposition \ref{lem:gl3}-$(ii)$, one finds that $$\Phi''(\bv_2)[\bh_0]^2 =\| V_2\|_2^2 -\b\intR V_2^3dx=(1-2\b)\| V_2\|_2^2<0\quad \mbox{provided}\quad \b>\frac 12.$$ See also Remark \ref{rem:2}.
\end{remark}
%%%%%%%%%%%%%%%%%%%%%%%%%%%%%%%%%%%%%%%%%%%%%%%%%%%%%%%%%%%%%%%%%%%%%%%%%%%
\section{Main results}\label{sec:mainR}
%%%%%%%%%%%%%%%%%%%%%%%%%%%%%%
\subsection{Existence of Ground states}\label{sec:ground}
%%%%%%%%%%%%%%%%%%%%%%%%%%%%%%
Concerning the existence of ground state solutions of
\eqref{NLS-KdV2}, the first result is the following.
\begin{Theorem}\label{th:1}
Suppose that $\b>\L$, then System \eqref{NLS-KdV2} has a positive
even ground state $\wt{\bu}=(\wt{u},\wt{v})$.
\end{Theorem}
\begin{pf}
We divide the proof into two steps. In the first step, we prove that $\inf_{\cN}\Phi$ is achieved at some positive function $\wt{\bu}\in\h$, while in the second step, we show that $\wt{\bu}$ is indeed a ground state, i.e.,
$$
\Phi(\wt{\bu})=\min\{ \Phi(\bu)\, :\: \bu\in\E,\:\Phi'(\bu)=0\}.
$$
{\it Step 1.}
By the Ekeland's variational principle; see \cite{eke}, there exists a PS sequence $\{\bu_k\}_{k\in \mathbb{N}}\subset\cN$, i.e.,
\be\label{eq:PS1}
\Phi (\bu_k)\to c=\inf_{\cN}\Phi
\ee
\be
\n_{\cN}\Phi(\bu_k)\to 0.
\ee
By \eqref{eq:restriction0}, easily one finds that $\{\bu_k\}$ is a bounded sequence on $\E$, and relabeling, we can assume that
$\bu_k\rightharpoonup \bu$ weakly in $\E$, $\bu_k\to \bu$ strongly in $\mathbb{L}^q_{loc}(\R)=L^q_{loc}(\R)\times  L^q_{loc}(\R)$ for every $1\le q<\infty$
and $\bu_k\to \bu$ a.e. in $\R^2$.
Moreover,  the constrained gradient  $\n_{\cN}\Phi (\bu_k)=\Phi' (\bu_k)-\eta_k \Psi'(\bu_k)\to 0$,  where $\eta_k$ is the corresponding Lagrange multiplier. Taking the scalar product with $\bu_k$ and recalling that  $(\Phi'(\bu_k)\mid \bu_k)=\Psi(\bu_k)=0$, we find that
$\eta_k (\Psi'(\bu_k)\mid \bu_k)\to 0$ and this
jointly with \eqref{eq:gamma}-\eqref{eq:bound} imply  that $\eta_k\to 0$. Since in addition $\|\Psi'(\bu_k)\|\leq C<+\infty$, we deduce that
$\Phi'(\bu_k)\to 0$.

Let us define $\mu_k=u_k^2+v_k^2$, where $\bu_k=(u_k,v_k)$. We {\it claim} that there is no evanescence, i.e., exist $R, C>0$ so that
\be\label{eq:vanishing}
\sup_{z\in\R}\int_{|z|<R}\mu_k\ge C>0,\quad\forall k\in\mathbb{N}.
\ee
On the contrary, if we suppose
$$
\sup_{z\in\R}\int_{|z|<R}\mu_k\to 0,
$$
by Lemma \ref{lem:measure}, applied in a similar way as in \cite{c-fract},  we find that $\bu_k\to \bo$ strongly in
$\mathbb{L}^{q}(\R)$ for any $2<q<\infty$, and as a consequence the weak limit $\bu^*\equiv \bo$. This is a contradiction since $\bu_k \in\cN$, and by \eqref{eq:restriction0}, \eqref{eq:restriction}, \eqref{eq:PS1} there holds
$$
0< \frac 17\rho <c+o_k(1)=\Phi(\bu_k)=F(\bu_k),\quad\mbox{with } o_k(1)\to 0 \quad\mbox{as }k\to\infty,
$$
hence \eqref{eq:vanishing} is true and the {\it claim} is proved.

We observe that   we can find a sequence of points
$\{z_k\}\subset\R^2$ so that by \eqref{eq:vanishing}, the translated sequence $\overline{\mu}_k(x)= \mu_k(x+z_k)$ satisfies
$$
\liminf_{k\to\infty}\int_{B_R(0)}\overline{\mu}_k\ge C >0.
$$
Taking into account that $\overline{\mu}_k\to \overline{\mu}$
strongly in $L_{loc}^1(\R)$, we obtain that
$\overline{\mu}\not\equiv 0$. Therefore, defining
$\overline{\bu}_k(x)=\bu_k(x+z_k)$, we have that  $\overline{\bu}_k$
is also a PS sequence of  $\Phi$ on $\cN$, in particular the weak
limit of $\overline{\bu}_k$, denoted by $\overline{\bu}$, is a
non-trivial critical point of $\Phi$ constrained on $\cN$, so
$\overline{\bu}\in\cN$. Thus, using \eqref{eq:restriction0} again,
we find
$$
\begin{array}{rcl}
\Phi (\overline{\bu}) & = & \dyle F(\overline{\bu})\\
& \le & \dyle\liminf_{k\to\infty} F(\overline{\bu}_k)\\
 & = & \dyle\liminf_{k\to\infty}\Phi(\overline{\bu}_k)= c.
 \end{array}
$$
Furthermore, by Proposition \ref{lem:gl3}-$(ii)$ we know that necessarily
$\Phi(\overline{\bu})<\Phi(\bv_2)$.

Taking into account that $\overline{\bu}\in\cN$, and the maximum principle, then $\overline{v}>0$, thus it is not difficult to show  $\wt{\bu}=|\overline{\bu}|=(|\overline{u}|,|\overline{v}|)=(|\overline{u}|,\overline{v})\in\cN$ with
\be\label{eq:min}
\Phi(\wt{\bu})=\Phi (\overline{\bu})=\min\{\Phi(\bu)\, :\: \bu\in\cN\},
\ee
so we have  $\wt{\bu}\ge \bo$. Finally, by the maximum principle applied to the first equation and the fact that $\Phi(\wt{\bu})<\Phi(\bv_2)$, we get $\wt{\bu}> \bo$.

\

{\it Step 2.} Assume, for a contradiction, there exists $\bw_0\in \E$  a non-trivial critical point of $\Phi$ such that
\begin{equation}\label{eq:sym}
\Phi(\bw_0)<\Phi(\wt{\bu})=\min\{ \Phi(\bu)\, :\: \bu\in\cN\}.
\end{equation}
Setting $\bw =|\bw_0|$ there holds
\begin{equation}\label{eq:sym1}
\Phi(\bw)=\Phi(\bw_0),\qquad {\Psi}(\bw)={\Psi}(\bw_0).
\end{equation}
For $\bw=(w_1,w_2)$, we set $\bw^\star=(w_1^\star,w_2^\star)$, where $w_j^\star$ is the Schwartz symmetric function associated to $w_j\ge 0$; $j=1,2
$.
Then by the classical properties of the Schwartz symmetrization; see for instance \cite{k},  there hold
\be\label{eq:primera}
\|\bw^{\star}\|^2\le \|\bw\|^2, \qquad G_\b(\bw^{\star})\ge G_\b(\bw),
\ee
thus, in particular, ${\Psi}(\bw^\star)\le {\Psi}(\bw)$.
Using the second identity of \eqref{eq:sym1} and the fact that $\bw_0$ is a critical point
of $\Phi$, we get ${\Psi}(\bw)={\Psi}(\bw_0)=0$. Furthermore,  there exists a unique $t_0>0$ so that $t_0\,\bw^{\star}\in {\mathcal{N}}$. In fact, $t_0$ comes from  $\Psi (t_0\bw^\star)=0$, i.e.,
\be\label{eq:t}
\| \bw^{\star}\|^2= t_0^2\intR (w_1^\star)^4dx+t_0\left( \frac 12\intR (w_2^\star)^3dx+\frac 32\b \intR (w_1^\star)^2w_2^\star \,dx\right),
\ee
then using that $\Psi(\bw)=0$, \eqref{eq:primera}-\eqref{eq:t}  and the fact that $\bw>\bo$ and $t_0>0$ we find
$$
\begin{array}{rcl}
& & \dyle \intR w_1^4\,dx+ \frac 12\intR w_2^3\,dx+\frac 32\b \intR w_1^2w_2\,dx \\ & & \\
& \ge &\dyle  t_0^2\intR w_1^4dx+t_0\left( \frac 12\intR w_2^3\,dx+\frac 32\b \intR w_1^2w_2 \,dx\right).
\end{array}
$$
Thus, clearly $t_0\le 1$, and as a consequence,
\be\label{eq:previa}
\Phi(t_0\,\bw^\star)= \frac 16 t_0^2\|\bw^\star\|^2+\frac{1}{12}t_0^4\intR (w_1^\star)^4\,dx\leq \frac 16 \|\bw\|^2+\frac{1}{12}\intR w_1^4\,dx=\Phi(\bw).
\ee
Therefore, inequalities \eqref{eq:previa}, \eqref{eq:sym} and the first identity of \eqref{eq:sym1}  yield
$$
\Phi(t_0\,\bw^{\star})\leq \Phi(\bw)<\Phi(\wt{\bu})= \min \{\Phi(\bu)\, :\:\bu\in \cN\},
$$
which is a contradiction because  $t_0\,\bw^{\star}\in \cN$.
\end{pf}
%%%%%%%%%%%%%%%%%%%%%%%%%%%%%%%%%%%%%%%%%%%%%%%%%%%%%%%%%%%%%%%%%%%%%%%%%%%%%%%%%
\begin{remark}\label{rem:2}
As we anticipated at the introduction, see also Remark \ref{rem:1},
in the setting by \cite{dfo}, $\l_1=\l_2$ and $\b>\frac 12$, we have
found positive even ground state solutions in contrast with the
non-negative bound states founded by \cite{dfo}.
\end{remark}
The last result in this subsection deals with the existence of
positive  ground states of  \eqref{NLS-KdV2} not only for $\b>\L$,
but also for $0<\b\le\L$, at least for  $\l_2$ large enough.
\begin{Theorem}\label{th:ground2}
There exists $\L_2>0$ such that if $\l_2>\L_2$, System \eqref{NLS-KdV2} has an even ground state
$\wt\bu>\bo$ for every $\b>0$.
\end{Theorem}
\begin{pf} Arguing in the same way as in the proof of Theorem \ref{th:1},
we initially have that  there exists an even ground state
$\wt{\bu}\ge \bo$. Moreover, in Theorem \ref{th:1} for $\b>\L$ we
proved that $\wt{\bu}>\bo$. Now we need to show that for $\b\le \L$
indeed $\wt{\bu}>\bo$ which follows by the maximum principle
provided $\wt{\bu}\neq \bv_2$. Taking into account Proposition
\eqref{lem:gl3}-$(i)$, $\bv_2$ is a strict local minimum, but this
does not allow us to prove that $\wt{\bu}\neq \bv_2$. The new  idea
here consists on proving the existence of  a function
$\bu_1=(u_1,v_1)\in\cN$ with $\Phi(\bu_1)<\Phi(\bv_2)$. To do so,
since $\bv_2=(0,V_2)$ is a local minimum of $\Phi $ on $\cN$
provided $0<\b<\L$, we cannot find $\bu_1$ in a neighborhood of
$\bv_2$ on $\cN$. Thus, we define $\bu_1=t(V_2,V_2)$ where  $t>0$ is
the unique value so that $\bu_1\in \cN$.

Notice that  $t>0$  is given by $\Psi(\bu_1)=0$, i.e.,
\be\label{eq:v22} \| (V_2,V_2)\|^2=t^2\intR V_2^4\, dx+\frac 12
t(1+3\b)\intR V_2^3\,dx.
\ee
Moreover,
\be\label{eq:v2}
\|(V_2,V_2)\|^2=2\|V_2\|_2^2+(\l_1-\l_2)\intR V_2^2\,dx=\intR V_2^3\,dx+(\l_1-\l_2)\intR V_2^2\,dx,
\ee
hence,
substituting \eqref{eq:v2} into \eqref{eq:v22} we get
$$
t^2\intR V_2^4\, dx+\frac 12 t(1+3\b)\intR V_2^3\,dx=\intR V_2^3\,dx+(\l_1-\l_2)\intR V_2^2\,dx,
$$
therefore, dividing the above expression by the $L^1$ norm of $V_2^3$, using that
$$
\intR \cosh^{-8}(x)\,dx=\frac{32}{35},\qquad \intR \cosh^{-6}(x)\,dx=\frac{16}{15},\qquad \intR \cosh^{-4}(x)\,dx=\frac{4}{3},
$$
and the definition of $V_2$ by \eqref{eq:segunda}, we find
\be\label{eq:segundo1}
\frac{18}{7}\l_2 t^2+\frac 12 t(1+3\b)-\left(1+5\frac{\l_1-\l_2}{12\l_2}\right)=0.
\ee
The energies  of $\bu_1$, $\bv_2$ are given by
$$
\Phi(t(V_2,V_2))=\frac 16 t^2\left(\intR V_2^3\,dx+ (\l_1-\l_2)\intR V_2^2\,dx\right)+\frac{1}{12}t^4\intR V_2^4\,dx,
$$
$$
\Phi (\bv_2)=\frac{1}{12}\intR V_2^3\,dx.
$$
Thus, we want to prove that for the unique $t>0$ given by \eqref{eq:v2} we have
$$
\frac 16 t^2\left(\intR V_2^3\,dx+ (\l_1-\l_2)\intR V_2^2\,dx\right)+\frac{1}{12}t^4\intR V_2^4\,dx<\frac{1}{12}\intR V_2^3\,dx,
$$
then arguing as for \eqref{eq:segundo1}, it is sufficient to prove that the following inequality holds
\be\label{eq:segundo2}
\frac{18}{7}\l_2 t^4+t^2\left(2+5\frac{\l_1-\l_2}{6\l_2}\right)-1<0.
\ee
Using \eqref{eq:segundo1} and the fact that  $2+5\frac{\l_1-\l_2}{6\l_2}>0$ for every $\l_1,\,\l_2>0$, fixed $\b>0$ we have that \eqref{eq:segundo2} is satisfied provided $\l_2$ is sufficiently large, namely $\l_2>\L_2>0$, proving that  $\Phi(\bu_1)<\Phi(\bv_2)$ which
 concludes the result.
 \end{pf}
%%%%%%%%%%%%%%%%%%%%%%%%%%%%%%
\subsection{Existence of Bound states}\label{sec:bound}
%%%%%%%%%%%%%%%%%%%%%%%%%%%%%%
In this subsection we establish  existence of bound states to
\eqref{NLS-KdV2}. The first theorem deals with  a perturbation
framework, in which we suppose that $\b=\e\wt{\b}$, with $\wt{\b}$
fixed and independent of $\e$. Note that $\wt{\b}$ can be negative, and $0<\e\ll1$. Then we rewrite the energy functional
$\Phi$ as $\Phi_\e$ to emphasize its dependence  on $\e$,
$$
\Phi_\e(\bu)=\Phi_0(\bu)-\tfrac 12\e \wt{\b} \intR u^2v\,dx,
$$
where $\Phi_0=I_1+I_2$.

Let us set $\bu_0=(U_1,V_2)$, where $V_2$ is given by
\eqref{eq:segunda} and $U_1$ is the unique positive solution of
$-u''+\l_1u=u^3$ in $H$; see \cite{cof,kw}. This function $U_1$ has
the following explicit expression, \be\label{eq:segunda2}
U_1(x)=\frac{\sqrt{2\l_1}}{\cosh (\sqrt{\l_1}x)}. \ee Note also that
$U_1$ satisfies the following identity, \be\label{eq:gr2}
 \|U_1\|_1=\inf_{u\in H\setminus\{0\}} \frac{\|u\|_1^2}{\left( \intR u^4dx\right)^{1/2}}.
\ee
\begin{Theorem}\label{th:2}
There exists $\e_0>0$ so that for any $0<\e <\e_0$ and
$\b=\e\wt{\b}$, System \eqref{NLS-KdV2} has an even bound state
$\bu_{\e}$ with  $\bu_{\e}\to \bu_0$ as $\e\to 0$. Moreover, if $\b>0$ then $\bu_\e>\bo.$
\end{Theorem}
In order to prove this result, we follow some ideas of \cite[Theorem
4.2]{c} with appropriate modifications.

\

\noindent {\it Proof of Theorem \ref{th:2}.}
It is well known that $U_1$ and $V_2$ are non-degenerate critical points of $I_1$ and $I_2$ on $H$ respectively; \cite{kw}. Plainly, $\bu_0$ is a non-degenerate critical point of $\Phi_0$ acting  on $\h$.
Then, by the Local Inversion Theorem, there exists a critical point $\bu_\e$ of $\Phi_\e$ for any $0<\e<\e_0$ with $\e_0$ sufficiently small; see  \cite{a-m} for more details. Moreover, $\bu_\e\to \bu_0$ on $\mathbb{H}$ as $\e\to 0$. To complete the proof it remains to show that if $\b>0$, then $\bu_\e> \bo$.

Let us  denote the positive part  $\bu_\e^+=(u_{\e}^+,v_{\e}^+)$ and  the negative part $\bu_\e^-=(u_{\e}^-,v_{\e}^-)$. By  \eqref{eq:gr2} we have
\begin{equation}\label{eq:pm1}
\|u_{\e}^\pm\|_1^2 \geq  \|U_1\|_1\left(\intR (u_{\e}^\pm)^4dx\right)^{1/2}.
\end{equation}

Multiplying  the second equation of \eqref{NLS-KdV2} by $v_{\e}^-$ and integrating on $\R$ one obtains
\be
\|v_{\e}^-\|_2^2 =\intR (v_{\e}^-)^3dx + \e \wt{\b}\intR (u_{\e})^2v_{\e}^- dx\le 0,
\ee
thus $\|v_\e^-\|_2=0$ which implies $v_\e=v_\e^+\ge 0$. Furthermore,  $\bu_\e\to\bu_0$ implies
$v_\e\to V_2$, which jointly with the maximum principle gives $v_\e>0$ provided $\e$ is sufficiently small.

Multiplying now  the first equation of \eqref{NLS-KdV2} by $u_{\e}^\pm$ and integrating on $\R$ one obtains
\begin{eqnarray*}
\|u_{\e}^\pm\|_1^2 &=& \intR (u_{\e}^\pm)^4dx + \e \wt{\b}\intR (u_{\e}^\pm)^2v_{\e} \,dx\\
&\leq&\intR (u_{\e}^\pm)^4dx +\e \wt{\b}\left( \intR (u_{\e}^\pm)^4dx\right)^{1/2}\left(   \intR v_{\e}^2 \,dx \right)^{1/2}.
\end{eqnarray*}
This, jointly with \eqref{eq:pm1}, yields
\be\label{eq:primera1}
\|u_{\e}^\pm\|_1^2 \le \frac{\|u_{\e}^\pm\|_1^4}{ \|U_1\|_1^2} + \e\, \theta_\e\;\frac{\|u_{\e}^\pm\|_1^2}{ \|U_1\|_1},
\ee
where
$$
\theta_\e = \wt{\b}\left( \intn v_{\e}^2\right)^{1/2}.
$$
Hence, if  $\|u_{1\e}^\pm\|>0$, one infers
\begin{equation}\label{eq:pmm1}
\|u_{\e}^\pm\|_1^2 \ge  \|U_1\|_1^2 +o(1),
\end{equation}
where $o(1)=o_\e (1)\to 0$ as $\e\to 0$.
Using again  $\bu_\e \to \bu_0$, then $u_{\e}\to U_1>0 $, as a consequence, for $\e$ small enough, $\|u_{\e}^+\|>0$. Thus
\eqref{eq:pmm1} gives
\begin{equation}\label{eq:pm11}
\|\bu_{\e}^+\|^2 =\|u_{\e}^+\|_1^2+\|v_{\e}^+\|_2^2 \geq   \|U_1\|_1^2 +o(1).
\end{equation}
Now, suppose for a contradiction,  that $\|u_{\e}^-\|_1>0$. Then as for \eqref{eq:pm11}, one obtains
\begin{equation}\label{eq:pm12}
\|\bu_\e^-\|^2 =\|u_{\e}^-\|_1^2+\|v_{\e}^-\|_2^2 \geq  \|U_1\|_1^2+o(1).
\end{equation}
On one hand, using \eqref{eq:pm11}-\eqref{eq:pm12}, we find
\be\label{eq:eval-funct}
\begin{array}{rcl}
\Phi(\bu_\e) & = & \dyle \tfrac 16 \|\bu_\e\|^2+\tfrac{1}{12}\intR u_\e^4\,dx \\ & &\\
& = &  \dyle\tfrac 16 \left[  \|\bu_\e^+\|^2 + \|\bu_\e^-\|^2\right]+\tfrac{1}{12}\intR [(u_\e^+)^4+(u_\e^-)^4]\,dx\\ & & \\
 &  \ge & \dyle \tfrac 16  \|\bu_0\|^2 +\tfrac 16 \| U_1\|_1^2+ \tfrac{1}{12}\intR U_1^4\,dx  +o(1).
\end{array}
\ee
On the other hand,  since $\bu_\e \to \bu_0$ we have
\be\label{eq:final}
\Phi(\bu_\e)= \tfrac 16 \|\bu_\e\|^2+\tfrac{1}{12}\intR u_\e^4\,dx\to \tfrac 16\|\bu_0\|^2+ \tfrac{1}{12}\intR U_1^4\,dx,
\ee
which is in contradiction with \eqref{eq:eval-funct}, proving that $u_{\e}\geq 0$.

In conclusion, we have proved that $v_\e>0$ and $u_\e\ge 0$. To
prove the positivity of $u_\e$, using once more that $\bu_\e\to
\bu_0$,  and $\b=\e\wt{\b}\ge 0$ we can apply the maximum principle
to the first equation of \eqref{NLS-KdV2}, which implies that
$u_\e>0$, and finally, $\bu_\e>\bo$. \rule{2mm}{2mm}\medskip

From  the existence of a positive ground state established in
Theorem \ref{th:1} for $\b>\L$, and more precisely in Theorem
\ref{th:ground2} for $\b>0$, provided  $\l_2$ is sufficiently large,
we can show the existence of a different positive bound state of
\eqref{NLS-KdV2} in the following.
\begin{Theorem}\label{th:bound2}
In the hypotheses of Theorem \ref{th:ground2} and $0<\b<\L$, there
exists an even bound state $\bu^*>\bo$ with
$\Phi(\bu^*)>\Phi(\bv_2)$.
\end{Theorem}
\begin{pf}
The positive ground state $\wt\bu$ founded in Theorem
\ref{th:ground2} satisfies $\Phi(\wt{\bu})<\Phi(\bv_2)$ and even
more, if $\b<\L$ by Proposition \ref{lem:gl3}, $\bv_2$ is a strict
local minimum of $\Phi$ constrained on $\cN$. As a consequence, we
have the Mountain Pass (MP in short) geometry between $\wt{\bu}$ and
$\bv_2$ on $\cN$. We define the set of all continuous paths joining
$\wt{\bu}$ and $\bv_2$  on the Nehari manifold by
$$
\G=\{ \g:[0,1]\to\cN\mbox{ continuous }|\: \g(0)=\wt{\bu},\: \g(1)=\bv_2\}.
$$
Thanks to the  MP Theorem by Ambrosetti-Rabinowitz; \cite{ar}, there exists a PS sequence $\{\bu_k\}\subset\cN$, i.e.,
$$
\Phi (\bu_k)\to c=\inf_{\cN}\Phi,\qquad \n_{\cN}\Phi(\bu_k)\to 0,
$$
where
\be\label{eq:MP-level}
c=\inf_{\g\in\G}\max_{0\le t\le 1}\Phi(\g(t)).
\ee
Plainly, by \eqref{eq:restriction0} the sequence $\{\bu_k\}$ is bounded  on $\h$, and we obtain a weakly convergent subsequence  $\bu_k\rightharpoonup\bu^*\in \cN$.

The  difficulty of the lack of compactness, due to work in the one dimensional case (see Remark \ref{rem:obs1}-$(ii)$), can be circumvent  in a similar way as in the proof of Theorem \ref{th:1}, so we omit the full detail for short. Thus, we find that the weak limit $\bu^*=(u^*,v^*)$ is an even bound state of \eqref{NLS-KdV2}, and clearly, $\Phi (\bu^*)>\Phi(\bv_2)$.

It remains to prove that  $\bu^*>\bo$. To do so, let us introduce the following problem
\begin{equation}\label{NLS-KdV+}
\left\{\begin{array}{rcl}
-u'' +\l_1 u & = & (u^+)^3+\beta u^+v \\
-v'' +\l_2 v & = & \frac 12 v^2+\frac 12\beta (u^+)^2.
\end{array}\right.
\end{equation}
By the maximum principle every nontrivial solution $\bu=(u,v)$ of
\eqref{NLS-KdV+}  has the second component  $v>0$ and the first one
$u\ge 0$.  Let us define its energy functional
$$
\Phi^+ (\bu)=\frac 12\|\bu \|^2 -G_\b (u^+,v),
$$
and consider the corresponding Nehari manifold
$$
\cN^+=\{\bu\in\h\setminus\{\bo\}\, :\: (\n\Phi^+ (\bu)|\bu)=0\}.
$$
Also, we denote
$$
I_1^+(u)=\frac12 \| u\|_1^2-\frac 14\intR (u^+)^4\,dx.
$$
It is not very difficult to show that the properties proved for
$\Phi$ and $\cN$ still hold for $\Phi^+$ and $\cN^+$. Unfortunately,
$\Phi^+$ is not $\mathcal{C}^2$, thus Proposition
\ref{lem:gl3}-$(i)$ does not holds directly for $\Phi^+$. To solve
this difficulty, we are going to prove that $\bv_2$ is a strict
local minimum of $\Phi^+$ constrained on $\cN^+$ without using the
second derivative of the functional. Note that in a similar way as
in \eqref{eq:tang1}, there holds \be \bh=(h_1,h_2)\in
T_{\bv_2}\cN^+\Longleftrightarrow h_2\in T_{V_2}\cN_2. \ee Taking
$\bh\in T_{\bv_2}\cN^+$ with $\|\bh\|=1$, we consider $\bv_\e=(\e
h_1,V_2+\e h_2)$. Plainly, there exists a unique $t_{\e}>0$ so that
$t_\e\bv_\e\in \cN^+$. Thus, we want to prove there exists $\e_1>0$
so that
$$
\Phi^+(t_\e\bv_\e)>\Phi^+ (\bv_2)\qquad \forall\: 0<\e<\e_1.
$$
It is convenient to distinguish if $h_1=0$ or not. In the former case, $h_1=0$, $\bv_\e=(0, V_2+\e h_2)$. Hence $t_\e\bv_\e\in \cN^+\Leftrightarrow t_\e(V_2+\e h_2)\in\cN_2$. Furthermore,
\be\label{eq:anterior}
\Phi^+(t_\e\bv_\e)=I_2(t_\e(V_2+\e h_2))>I_2(V_2)=\Phi(\bv_2)=\Phi^+(\bv_2),
\ee
where the previous inequality holds because $V_2$ is a strict local minimum of $I_2$  on $\cN_2$.

Let us now consider the case $h_1\neq 0$. There holds
\be\label{eq:posterior}
\Phi^+(t_\e \bv_\e)=I_2(t_\e(V_2+\e h_2))+I_1^+(t_\e\e h_1)-\frac 12\b \e^2t_\e^2\intR (h_1^+)^2(V_2+\e h_2)\,dx.
\ee
By \eqref{eq:anterior} and \eqref{eq:posterior} it follows,
\be\label{eq:posterior2}
\Phi^+(t_\e \bv_\e)>\Phi^+(\bv_2)+I_1^+(t_\e\e h_1)-\frac 12\b \e^2t_\e^2\intR (h_1^+)^2(V_2+\e h_2)\,dx.
\ee
To finish, it is sufficient to show that
$$
\mathcal{J} (t_\e\bv_\e):=I_1^+(t_\e\e h_1)-\frac 12\b \e^2t_\e^2\intR ( h_1^+)^2(V_2+\e h_2)\,dx>0\qquad \forall\: 0<\e<\e_1.
$$
Let $\a<1$ be such that $\a>\frac{\b}{\L}$. By \eqref{eq:Lambda} and $\b<\L$ there holds
$$
\b\intR V_2( h_1^+)^2\,dx <\a\| h_1\|_1^2,
$$
then for $\e_1$ smaller than before (if necessary) we have
\be\label{eq:posterior3}
\b\intR (V_2+\e h_2)(h_1^+)^2\,dx <\a\| h_1\|_1^2 \qquad \forall\: 0<\e<\e_1.
\ee
Using  \eqref{eq:posterior3} and the Sobolev inequality, we obtain
$$
\mathcal{J} (t_\e\bv_\e)>\frac 12t_\e^2\e^2\| h_1\|_1^2(1-t_\e\a-ct_\e^2\e^2),\quad\mbox{ for a constant } c>0.
$$
Now, taking into account that $t_\e\to 1$ as $\e\searrow 0$, we infer there exists a  constant $c_0>0$ so that
\be\label{eq:final2}
\mathcal{J} (t_\e\bv_\e)>\e^2c_0\|h_1\|_1^2.
\ee
Finally, by \eqref{eq:posterior2}, \eqref{eq:final2} it follows that
$$
\Phi^+(t_\e\bv_\e)>\e^2c_0\| h_1\|_1^2+\Phi^+(\bv_2) > \Phi^+(\bv_2),
$$
which proves that $\bv_2$ is a strict local minimum for $\Phi^+$ on $\cN^+$.

\noindent From the preceding arguments, it follows that $\Phi^+$ has a MP critical point $\bu^*\in \cN^+$, which gives rise to a solution of \eqref{NLS-KdV+}.
In particular, one finds that $u, v\ge 0$.
   In addition, since $\bu^*$ is a MP critical point, one
has that $\Phi (\bu^*)=\Phi^+(\bu^*)>\Phi^+(\bv_2)=\Phi(\bv_2)>0$,
which implies  $u^*\ge 0$ with $u^*\not\equiv 0$, and by the maximum
principle applied to each single equation we get $u^*,\, v^*>0$,
hence $\bu^*>\bo$. \end{pf}

In view of Theorems \ref{th:ground2}, \ref{th:bound2}, some remarks
are in order.
\begin{Remarks}\label{rem:11}
\begin{itemize}
\item[(i)] Following the proof of Theorem \ref{th:bound2},  a natural question is what happens in the limit
case $\b=\L$. In that case $\bu^*$ could coincides with $\bv_2$ which is
non-negative, but not positive. Indeed, this is our conjecture in
view of the second equation by \eqref{eq:lim}; see also Figure
 \ref{figure1}.
\item[(ii)] In the hypotheses of Theorems \ref{th:ground2}, \ref{th:bound2} we have found the coexistence of two positive solutions, the ground state $\wt{\bu}$ in Theorem \ref{th:ground2} and the bound state $\bu^*$ in Theorem \ref{th:bound2}, proving a non-uniqueness result of positive solutions to \eqref{NLS-KdV2}. This is a great difference with the more studied system of coupled NLS equations
$$
\left \{
\begin{array}{ll}
- \D u_1+ \l_1 u_1 &=  \mu_1u_1^3+\b u_2^2 u_1\\
- \D u_2 + \l_2 u_2 &= \mu_2  u_2^3+\b u_1^2u_2,
\end{array} \right.
$$
(see for instance
\cite{a,ac1,ac2,bt,chen-zou,c,fl,iko,itanaka,linwei,linwei2,liu-wang,mmp,sirakov,ww,wy}
and the references therein) for which it is known that there is
uniqueness of positive solutions, under appropriate conditions on
the parameters including the case $\b>0$ small; see more
specifically  \cite{iko,wy}. Indeed, for $\b>0$ small, the ground
state is not positive, and it is given by one of the two
semi-trivial solutions  $(U^{(1)},0)$ or $(0,U^{(2)})$ depending on
if $\Phi (U^{(1)},0)$ is lower or grater than $\Phi(0,U^{(2)})$
which plainly corresponds to $\l_1^{2-\frac n2}\mu_2<\l_2^{2-\frac
n2}\mu_1$ or $\l_1^{2-\frac n2}\mu_2>\l_2^{2-\frac n2}\mu_1$
respectively. Here $U^{(j)}$ is the unique\footnote{See
\cite{cof,kw} for this uniqueness result.} positive radial solution
of $-\D u_j+\l_ju_j=\m_j u_j^3$ in $W^{1,2}(\mathbb{R}^n)$, for
$n=1,\,2,\, 3$ and $j=1,\, 2$.
\end{itemize}
\end{Remarks}
%%%%%%%%%%%%%%%%%%%%%%%%%%%%%%%%%%%%%%%%%%%%%%%%%%%%%%%%%%%%%%%%%%%%%%%%%%%%%%%%%
\section{An extended  NLS-KdV system with general power nonlinearities}\label{sec:extended}
%%%%%%%%%%%%%%%%%%%%%%%%%%%%%%%%%%%%%%%%%%%%%%%%%%%%%%%%%%%%%%%%%%%%%%%%%%%%%%%%%
In this section we want to show that if one consider a more general
system than \eqref{NLS-KdV}, \eqref{NLS-KdV2} with more general
power nonlinearities, like the following
\begin{equation}\label{NLS-KdV-pq-general}
\left\{\begin{array}{rcl}
if_t + f_{xx} + \tau_1\,|f|^{q-1}f + \b fg & = & 0\\
g_t +g_{xxx} +\tau_2 \,|g|^{p-1}g_x +\frac 12\b (|f|^2)_x  & = & 0,
\end{array}\right.
\end{equation}
where $\tau_1,\,\tau_2,\, \b$ are real constants, the one can prove
the same results of the previous section with appropriate
hypotheses. Looking for solutions of \eqref{NLS-KdV-pq-general} in
the form by \eqref{eq:sol-t-v-solutions} we find that for
$\mu_1=\tau_1$, $\mu_2=\frac{\tau_2}{p}$, the real functions $u,v$
solve the following system
\begin{equation}\label{NLS-KdV-pq-general2}
\left\{\begin{array}{rcl}
-u'' + \l_1 u & = & \mu_1|u|^{q-1}u +\b uv\\
-v'' +\l_2 v & = & \mu_2 |v|^{p-1}v+\frac 12\b u^2,
\end{array}\right.
\end{equation}
where we consider $\l_j,\, \m_j>0$; $j=1,2$; $p,\,q\ge 2$. We take $\b>0$ in order to obtain positive solutions, although some  results about the existence of bound states hold true too without the positivity of them.

Some of the results in Section \ref{sec:mainR} hold with minor changes. Notice that, since we look for positive solutions of \eqref{NLS-KdV-pq-general2}, one could consider the term $(|g|^p)_x$ (as in previous sections where $p=2$) instead of $|g|^{p-1}g_x$ in \eqref{NLS-KdV-pq-general}, and hence one would have $|g|^p$ instead of $|g|^{p-1}g$ in \eqref{NLS-KdV-pq-general2}, obtaining the same existence of positive bound and ground states that we will prove here,  in Theorems \ref{th:14}, \ref{th:ground-bound2-p-q}.  More general systems than  \eqref{NLS-KdV-pq-general2} will be analyzed in a forthcoming paper.

Note that \eqref{NLS-KdV-pq-general2} has a unique non-negative semi-trivial solution defined $\bv_{p}=(0,V_p)$ with $V_p$ the unique positive solution of $-v''+\l_2v=\m_2|v|^{p-1}v$ in $H$,  which have the following explicit expression,
\be\label{eq:v-p}
V_p(x)=\left[ \frac{(p+1)\l_2}{2\m_2\cosh^2\left(\frac{p-1}{2}\sqrt{\l_2}\, x\right)}\right]^{\frac{1}{p-1}}.
\ee
Following similar notation as for \eqref{NLS-KdV2}, we denote the associated energy functional of \eqref{NLS-KdV-pq-general2}  by
\be\label{eq:funct-p-q}
\Phi(\bu)=J_1(u)+J_2(v)- \tfrac 12\b \intR u^2v\,dx,\qquad \bu\in \E,
\ee
with
$$
J_1(u)=\tfrac 12 \|u\|_1^2 -\frac{\m_1}{q+1}\, \intR |u|^{q+1}dx,\quad J_2(v)=\tfrac 12 \|v\|_2^2 -\frac{\m_2}{p+1}\, \intR |v|^{p+1}dx;\qquad u,\, v\in E.
$$
Also, for $\Psi (\bu)=(\n \Phi(\bu)|\bu)$, we define the corresponding Nehari manifold as
$$
\cN =\{ \bu\in \h\setminus\{\bo\}: \Psi(\bu)=0\}.
$$
Plainly,
$$
(\n\Psi(\bu)|\bu)=-\|\bu\|^2-\m_1(q-2)\intR |u|^{q+1}\, dx-\m_2(p-2)\intR |v|^{p+1}\, dx\qquad\forall\, \bu\in \cN,
$$
thus $\cN$ is a smooth manifold locally near any point $\bu\not=
\bo$ with $\Psi(\bu)=0$. Moreover, $\Phi''(\bo)= I_1''(0)+I_2''(0)$
is positive definite, so we infer that $\bo$ is a strict minimum of
$\Phi$. As a consequence,  $\bo$ is an isolated point of the set
$\{\Psi(\bu)=0\}$, proving that, on one hand $\cN$ is a smooth
complete manifold of codimension $1$, and on the other hand there
exists a constant $\rho>0$ so that \be\label{eq:boundd}
\|\bu\|^2>\rho,\qquad\forall\,\bu\in \cN. \ee Then, as for
\eqref{NLS-KdV2} where  $q=3$, $p=2$ one has that $\bu\in
\h\setminus\{\bo\}$ is a critical point of $\Phi$ if and only if
$\bu\in\cN$ is a critical point of $\Phi$ constrained on $\cN$.
Furthermore,
$$
\Phi (\bu)=\frac 16\|\bu\|^2+\left(\frac 13-\frac{1}{q+1}\right)\m_1\intR|u|^{q+1}\, dx+\left(\frac 13-\frac{1}{p+1}\right)\m_2\intR|v|^{p+1}\, dx\qquad\forall\, \bu\in \cN,
$$
then clearly by \eqref{eq:boundd} and the previous identity, $\Phi$ on $\cN$ is bounded bellow, for every $2\le p<\infty$, $2\le q<\infty$.
\begin{Proposition}\label{prop:gl}
For  $\L$  defined by \eqref{eq:Lambda}:
\begin{itemize}
\item[(i)]  if $\b\le \L$, then $\bv_p$ is a strict local minimum  of $\Phi$ constrained on $\cN$,
\item[(ii)] for any   $\b>\L$, then $\bv_p$ is a  saddle point of $\Phi$ constrained on $\cN$. Even more,  $\dyle\inf_{\cN}\Phi <\Phi(\bv_2)$.
\end{itemize}
\end{Proposition}
The proof is a straightforward  calculation of the one of Proposition \ref{lem:gl3}. Furthermore, defining $U_q$ as the unique positive solution of $-u''+\l_1u=\m_1 |u|^{q-1}u$ in $H$ (given by \eqref{eq:v-p} substituting $p$ by $q$), we have the following.
\begin{Theorem}\label{th:14} Assume that $2\le p<\infty$, $2\le q<\infty$.
\begin{itemize}
\item[(i)] If $\b>\L$, then System \eqref{NLS-KdV-pq-general2} has a positive even ground state $\wt{\bu}=(\wt{u},\wt{v})$.
\item[(ii)] There exists $\e_0>0$ such that for any $0<\e <\e_0$ and $\b=\e\wt{\b}>0$, System \eqref{NLS-KdV-pq-general2} has an even bound state $\bu_{\e}> \bo$ with  $\bu_{\e}\to \bu_0=(U_q,V_p)$ as $\e\to 0$.
\end{itemize}
\end{Theorem}
\begin{pf}
We can adapt, with appropriate modifications, the ideas in the proof of Theorem \ref{th:1}, since the nonlinearity $|v|^{p+1}$ is even while in Theorem \ref{th:1} the nonlinearity on $v$ is $v^3$ (odd). That proves part $(i)$.

 Part $(ii)$ follows by a little modification of the ideas of Theorem \ref{th:2}, by the same reason as above.
\end{pf}

Concerning the existence of ground states for any $\b>0$ one can follow the proof of  Theorem \ref{th:ground2} that it holds true too with a restriction on the power exponent $q>2p-2$, which trivially holds for \eqref{NLS-KdV2} where $q=3,\,p=2$. Thus, using that property, it is not difficult to show also the existence of positive bound states for $0<\b<\L$ as in Theorem \ref{th:bound2}. We enunciate these results in the following.
\begin{Theorem}\label{th:ground-bound2-p-q}
Assume that $2\le p<\infty$, $2\le q<\infty$ and even more $q>2p-2$, then:
\begin{itemize}
\item[(i)] there exists $M>0$ such that if $\l_2>M$, System \eqref{NLS-KdV-pq-general2} has an even ground state $\wt\bu>\bo$ for every $\b>0$,
\item[(ii)] if $\l_2>M$ and $0<\b<\L$, there exists an even bound state $\bu^*>\bo$ with $\Phi(\bu^*)>\Phi(\bv_2)$.
\end{itemize}
\end{Theorem}
\begin{Remarks}
\begin{itemize}
\item[(i)] The restriction $q>2p-2$ appears when one tries to prove that $$\Phi(t(V_p,V_p))<\Phi (\bv_p)$$ for $t(V_p,V_p)\in\cN$. It does not seem to be optimal. Another test function different from $t(V_p,V_p)$ could circumvent this difficulty.
\item[(ii)] When $p=2$, $\m_2=p+1$, in  \cite[Theorem 4.1]{dfo2} Dias et al. impose $\b>3$ to obtain even bound states.
In our Theorem \ref{th:14}-$(i)$, following the idea of Remark
\ref{rem:1}, it is easy to see that it holds for  $\b> 3-a$ for some
constant $a>0$ when $\l_2>\l_1$, obtaining positive even bound and
ground states.
\item[(iii)] Note that here, in Theorems \ref{th:14}, \ref{th:ground-bound2-p-q} we have $2\le p<\infty$,
$2\le q<\infty$ obtaining positive even bound and ground states, in
contrast with \cite[Theorem 4.1]{dfo2} $2<q<5$, $p\in\{ 2,3,4\}$ and
$\m_2=p+1$ and  in \cite[Theorem 1.1]{ab} $2\le q<5$, $2\le p<5$
with $p$ a rational number with odd denominator, where the authors
obtained non-negative  even bound states in the former and positive
even bound states in the later.
\end{itemize}
\end{Remarks}
%%%%%%%%%%%%%%%%%%%%%%%%%%%%%%%%%
\section{Further results}\label{sec:remarks}
%%%%%%%%%%%%%%%%%%%%%%%%%%%%%%%%%
In this last section we show some results for explicit solutions. We
point out some remarks and open problems. To finish, we study some
extended systems with three or more equations.
%%%%%%%%%%%%%%%%%%%%%%%%%%%%%%%%%
%%%%%%%%%%%%%%%%%%%%%%%%%%%%%%%%%
\subsection{Explicit solutions}
In the particular case $0<\b<\frac 16$, $\l_2=4\l_1+\frac{1}{12}\b
(1-6\b)$ there exists a nontrivial  explicit solution\footnote{although the  results in this subsection can be established for the more general system \eqref{NLS-KdV-pq-general2}, we restrict ourselves to \eqref{NLS-KdV2} for short.}
$\bu_\b=(u_\b,v_\b)$ of \eqref{NLS-KdV2} defined by
$$
u_\b(x)=\frac{\sqrt{2\l_1(1-6\b)}}{\cosh(\sqrt{\l_1}x)},\qquad v_\b(x)=\frac{12\l_1}{\cosh^2(\sqrt{\l_1}x)}.
$$
Clearly one has that
\be\label{eq:lim}
\lim_{\b\searrow 0}\bu_\b= \bu_0=(U_1,V_2),\qquad \lim_{\b\nearrow \frac 16}\bu_\b=\bv_2=(0,V_2),
\ee
where $U_1$, $V_2$ are defined by  \eqref{eq:segunda2}, \eqref{eq:segunda} respectively. Then, the family $\{\bu_\b:\, 0<\b<\frac 16\}$ joins $\bu_0=(U_1,V_2)$ with $\bv_2$.
\begin{Remarks}
\begin{itemize}
\item[(i)]
If  we would have $\Phi(\bv_2)\ge\Phi(\bu_\b)$ in the range
$0<\b<\min\{ \frac 16, \L\}$, then we would be able to prove the
existence of a positive even bound state $\bu^*$ with
$\Phi(\bu^*)>\max\{\Phi(\bu_\b),\Phi(\bv_2)\}=\Phi(\bv_2)$, and in
particular we would have a non-uniqueness of positive solutions
result by a different way as in the previous sections.
Unfortunately, if $0<\b< \frac 16$ then \be\label{eq:v-segunda}
\begin{array}{rcl}
\Phi (\bv_2) & = & \dyle\frac 16\|\bv_2\|^2= \frac{1}{12}\intR V_2^3\,
dx=\frac{9}{2}\l_2^3\intR \frac{1}{\cosh^6(\frac{\sqrt{\l_2}}{2}x)}\,dx=\frac{24}{5}\l_2^{5/2}\\
\\
& = &\dfrac{24}{5}\left[4\l_1+\frac{1}{12}\b (1-6\b)\right]^{5/2},
\end{array}
\ee
and
\be\label{eq:beta}
\begin{array}{rcl}
\Phi (\bu_\b) & = & \dyle\frac 16\| \bu_\b\|^2+\frac{1}{12}\intR u_\b^4\,dx=\frac 16\left(\|u_\b\|_1^2+\|v_\b\|^2_2\right)+\frac{1}{12}\intR u_\b^4\,dx\\ & & \\
& = &\dfrac 16\left(\dyle\intR u_\b^4\,dx +\b \intR u_\b^2v_\b\, dx+\frac 12\left(\intR v_\b^3\,dx+\b\intR u_\b^2v_\b\, dx\right)\right)+\dyle\frac{1}{12}\intR u_\b^4\,dx\\ & & \\
& = &\dyle \frac 14\intR u_\b^4\,dx+\dfrac{\b}{4}\intR u_\b^2v_\b\,dx+\dfrac{1}{12}\intR v_\b^3\,dx\\ &  & \\
& = & \l_1^{2}(1-6\b)\dyle\intR\dfrac{1}{\cosh^4(\sqrt{\l_1}x)}\,dx+144\l_1^3\intR \dfrac{1}{\cosh^6(\sqrt{\l_1}x)}\,dx\\ &  &\\
& = & \dyle \frac{4}{3}\l_1^{3/2}(1-6\b)+\dfrac{768}{5}\l_1^{5/2}.
\end{array}
\ee
Comparing both energies, it is not difficult to show that $\Phi(\bv_2)<\Phi(\bu_\b)$ for every $0<\b<\frac 16$.
\item[(ii)] In  the setting in which the explicit solutions $\bu_\b$ of \eqref{NLS-KdV2} exists, i.e.,
$0<\b<\frac 16$, $\l_2=4\l_1+\frac{1}{12}\b (1-6\b)$, under
hypotheses of  Theorems  \ref{th:ground2} and \ref{th:bound2},  we have the  multiplicity of positive solutions: the ground
state $\wt\bu$ and the bound state $\bu^*$. We  conjecture that
$\bu_\b$ coincides with $\bu^*$, and hence, the suggestive
bifurcation diagram by Figure 1 holds not only for $\bu_\b$ but also for
$\bu^*$.
\begin{figure}[!htp]\label{figure1}
\hspace*{2cm}
\begin{center}
\includegraphics[width=320pt]{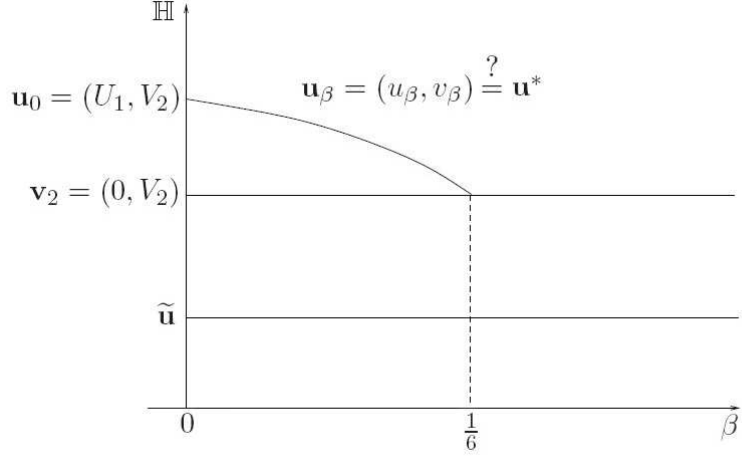}
\end{center}
\caption{Suggestive bifurcation diagram of the family
$\bu_\b=(u_\b,v_\b)$.}
\end{figure}
\end{itemize}
\end{Remarks}
\subsection{Some systems with more than two equations}
In this last subsection, we deal with some extended systems of \eqref{NLS-KdV2} to more than two\footnote{similar extensions of \eqref{NLS-KdV-pq-general2} to a more dimensional case, as \eqref{NLS-KdV2-n} extends \eqref{NLS-KdV2}, can be considered (at least in the subcritical framework with $ p,\,q<2^*$ defined in Remark \ref{rem:compacidad}), proving similar results as Theorems \ref{th:4}, \ref{th:5}.} equations, but also consider \eqref{NLS-KdV2} in a different dimensional case.

Note that  System \eqref{NLS-KdV}  has no sense in the dimensional case $n=2,\, 3$, however, \eqref{NLS-KdV2}  makes sense to be extended to more dimensions. Moreover, the results of the previous Sections can be established in the dimensional case $n=2,3$ with minor changes for system
\begin{equation}\label{NLS-KdV2-n}
\left\{\begin{array}{rcl}
-\D u +\l_1 u & = & u^3+\beta uv \\
-\D v +\l_2 v & = & \frac 12 v^2+\frac 12\beta u^2,
\end{array}\right.
\end{equation}
working on the corresponding Sobolev Spaces $E=W^{1,2}(\R^n)$, $n=2,3$ and its radial subspace $H=E_r$. In particular, Theorems \ref{th:1}, \ref{th:ground2}, \ref{th:2} and \ref{th:bound2} hold, obtaining the corresponding
 positive bound and ground state solutions which are radially symmetric in this case.

 \begin{Remarks}\label{rem:compacidad}
 \begin{itemize}
\item[(i)]
For $n=2,3$ there is no lack of compactness,
 because there holds the compact embedding of the radial Sobolev Space $H$ for all $2<s<2^*$ (see \cite{Lions-JFA82}), where $2^*=\infty$ if $n=2$ and $2^*=\frac{2n}{n-2}$ for $n=3$, which allow us to prove the Palais-Smale condition\footnote{in a similar way as in \cite[Lemma 3.2]{ac2}.} working on  $\h$.
\item[(ii)] Following some ideas by Ambrosetti and Colorado in \cite{ac2}, Liu and Zheng proved in \cite{lz} a partial result on existence of solutions to \eqref{NLS-KdV2-n} in the dimensional case $n=2,\, 3$. More precisely, in \cite{lz} the authors proved that the infimum of the functional associated to \eqref{NLS-KdV2-n} on the corresponding Nehari manifold is achieved, but they do not proved that it is positive, and it was not shown that the infimum on the Nehari Manifold is a ground state, i.e., the least energy solution of the functional as we have proved here for $n=1,\,2,\,3$. Also, in \cite{lz} was not investigated the existence of other bound states, as he have done in this manuscript and not only in the non-critical dimensions $n=2,\,3$ but also in the one dimensional case, $n=1$.
 \end{itemize}
\end{Remarks}
System \eqref{NLS-KdV2-n} can be seen as the stationary system of two coupled
NLS-NLS equations when one looks for solitary wave solutions, and  $(u,v)$ are the corresponding standing wave solutions. It is well known that systems of
NLS-NLS time-dependent equations have applications in some aspects of Optics, Hartree-Fock theory for Bose-Einstein
 condensates, among other physical phenomena; see for instance the earlier mathematical works
 \cite{Ack,a,ac1,ac2,acr,bt,fo,linwei,linwei2,mmp,sirakov,ww}, the more recent list (far from complete)
 \cite{chen-zou,itanaka,liu-wang} and references therein.

By the above discussion, one can motivate, from the application point of view, the study of the
following system of NLS-KdV-KdV equations,
\begin{equation}\label{NLS-KdV2-3}
\left\{\begin{array}{rcl}
-\D u +\l_0 u & = & u^3+\b_1 uv_1+\b_{2}uv_2 \\
-\D v_1 +\l_1 v_1 & = & \frac 12 v_1^2+\frac 12\beta_{1} u^2\\
-\D v_2+\l_2 v_2 & = & \frac 12v_2^2+\frac 12\beta_{2} u^2.
\end{array}\right.
\end{equation}
This system can also be seen as a perturbation of \eqref{NLS-KdV2-n} if $n=2,\,3$ or a perturbation of \eqref{NLS-KdV2} if $n=1$, when  $|\b_{1}|$ or $ |\b_{2}|$ is small.

Now, we use the same notation as in previous sections  with natural
meaning, for example, $\h=H\times H\times H$, $\E= E\times E \times
E$, $\bo=(0,0,0)$, \be\label{eq:Phi3} \Phi (\bu)=\frac 12 \|
\bu\|^2-\frac 14\intR u^4\, dx-\frac 16\intR (v_1^3+v_2^3)\,dx-\frac
12\intR u^2(\b_1v_1+\b_2v_2)\,dx \ee \be\label{eq:N3} \cN=\{\bu\in
\E\setminus\{\bo\}:  (\Phi'(\bu)| \bu)=0\}, \ee etc.

Let $U^*,\, V_j^*$ be the unique positive radial solutions of $-\D u+\l u=u^3$, $-\D v+\l_j v=\frac 12 v^2$
in $E$ respectively, $j=1,\,2$; see \cite{cof,kw}.
Then we have the following.
\begin{Theorem}\label{th:4}
There exists $\e_0>0$ such that for any $0<\e <\e_0$ and
$\b_{j}=\e\wt{\b}_{j}>0$, $j=1,2$, System \eqref{NLS-KdV2-3} has a
radial bound state $\bu_{\e}^*$ with  $\bu_{\e}\to
\bu_0^*=(U^*,V_1^*,V_2^*)$ as $\e\to 0$. Moreover, if $\b_j>0$ for $j=1,2$ then $\bu_{\e}^*>\bo$.
\end{Theorem}
The proof  follows  in a similar way as the proof of Theorem \ref{th:2}  with appropriate modifications, so that we omit it for short.

We also can prove the existence of a positive and radial ground
state of \eqref{NLS-KdV2-3} when the coupling parameters $\b_{j}$,
$j=1,2$ are sufficiently large. To do so, we define
\be\label{eq:Lambda-j} \L_j=\inf_{\varphi\in H\setminus\{
0\}}\frac{\|\varphi\|_0^2}{\intn V_{j}^*\varphi^2dx}\qquad j=1,2.
\ee where $\|\cdot\|_0$ is the norm in $E$ with $\l_0$.
\begin{Remark}\label{rem:semitrivial3}
The unique non-negative semi-trivial solutions of \eqref{NLS-KdV2-3} are given by
$\bv_1^*=(0,V_1^*,0)$, $\bv_2^*=(0,0,V^*_2)$ and $\bv_{12}^*=(0,V_1^*,V^*_2)$ with
$V_j^*$ for $j=1,2$ are defined before Theorem \ref{eq:Phi3}.
\end{Remark}
Concerning the ground states of  \eqref{NLS-KdV2-3}, the first result is the following.
\begin{Theorem}\label{th:5}
If $\b_{j}>\L_j$ for $j=1,\, 2$, then \eqref{NLS-KdV2-3} has a
positive radial ground state $\wt{\bu}$.
\end{Theorem}
\begin{pf}
As in Proposition \ref{lem:gl3}, using now that $\b_{j}>\L_j$,
j={1,\,2}, one can show that both $\bv_1^*$, $\bv_2^*$ are saddle
points of the energy functional $\Phi$  (defined by \eqref{eq:Phi3})
constrained on the Nehari manifold $\cN$ (defined by
\eqref{eq:Phi3}). Then \be\label{eq:energy3} c=\inf_{\cN}\Phi<\min
\{ \Phi(\bv_1^*),\Phi
(\bv_2^*)\}<\Phi(\bv_{12}^*)=\Phi(\bv_1^*)+\Phi (\bv_2^*). \ee By
the Ekeland's variational principle, there exists a PS sequence
$\{\bu_k\}_{k\in \mathbb{N}}\subset\cN$, i.e., \be\label{eq:PS1-3}
\Phi (\bu_k)\to c \ee \be \n_{\cN}\Phi(\bu_k)\to 0. \ee If $n=2,3$,
and $\{\bu_k\}$ satisfies the  PS condition; see Remark
\ref{rem:compacidad}-$(i)$. Thus, there exists a convergent subsequence
(denoted equal for short)  $\bu_k\to \wt{\bu}$. Arguing in a similar
way as in the proof of Theorem \ref{th:1} we have that $\wt{\bu}\ge
\bo$ and moreover, by the Schwartz symmetrization  properties,  one
can prove that
$$
c=\Phi(\wt{\bu})=\min\{ \Phi(\bu)\, :\: \bu\in\E,\:\Phi'(\bu)=0\}.
$$
To prove the positivity of $\wt{\bu}$, if one supposes that the
first component $u^* \equiv 0$, since the only non-negative
solutions of \eqref{NLS-KdV2-3} are the semi-trivial solutions
defined  in Remark \ref{rem:semitrivial3}, we obtain a contradiction
with \eqref{eq:energy3}. Furthermore, if the second or third
component vanish then $\wt{\bu}$ must be $\bo$, and this is not
possible because $\Phi|_{\cN}$ is bounded bellow by a positive
constant like in \eqref{eq:restriction}, then  $\bo$ is an isolated
point of the set $\{\bu\in \h\, :\: (\Phi'(\bu)|\bu)=0\}$, proving
that $\cN$ is a complete manifold. Finally, the maximum principle
gives that $\wt{\bu}>\bo.$

In the one dimensional case, $(n=1)$, the lack of compactness discussed in Remarks \ref{rem:obs1}-$(ii)$
can be circumvent  as in the proof of
Theorem \ref{th:1}, proving the result.
\end{pf}

\

Even more, one can show similar results to Theorems \ref{th:ground2}, \ref{th:bound2} in the setting of
 \eqref{NLS-KdV2-3}. To do so, we first prove an auxiliary result in the following.
\begin{Proposition}\label{prop:2}
Assume that $\b_{j}<\L_j$, $j=1,\, 2$ and moreover
$\dfrac{\b_{1}}{\L_{1}}+\dfrac{\b_{2}}{\L_{2}}<1$. Then $\bv_{12}^*$
is a strict local minimum of $\Phi$ constrained on $\cN$.
\end{Proposition}
\begin{pf}
Notice that since $\bv^*_{12}$ is a critical point of $\Phi$ then $D^2\Phi_{\cN}(\bv^*_{12})[\bh]^2=\Phi''(\bv_{12}^*)[\bh]^2$ for all $\bh\in T_{\bv_{12}^*}\cN$. We denote
$$
I_j(v)=\frac 12 \| v\|_{\l_j}^2-\frac 16\intn v^3dx.
$$
Then, using that $V_j^*$ is a strict local minimum of $I_j$ we have
there exist two positive constants $c_1,\,c_2$ so that for
$\bh=(h_0,h_1,h_2)\in T_{\bv_{12}^*}\cN$, \be\label{eq:Phi-segunda2}
\Phi''(\bv_{12}^*)[\bh]^2 \ge c_1\| h_1\|_1^2+c_2\| h_2\|_2^2+ \|
h_0\|^2-\b_{1}\intn h_0^2V_1^*\,dx-\b_{2}\intn h_0^2V_2^*\,dx. \ee
From \eqref{eq:Phi-segunda2} and
$\dfrac{\b_{1}}{\L_{1}}+\dfrac{\b_{2}}{\L_{2}}<1$ we infer that
$$
\Phi''(\bv_{12}^*)[\bh]^2 \ge c_1\| h_1\|_1^2+c_2\| h_2\|_2^2+
\left(1- \dfrac{\b_{1}}{\L_{1}}+\dfrac{\b_{2}}{\L_{2}}\right)\|
h_0\|^2,
$$
which proves the result.
\end{pf}
\begin{Theorem} Assume that $\b_{1},\b_{2}>0$, for $\l_1,\l_2$ large
enough:
\begin{itemize}
\item[(i)]  there exists a radial ground state $\wt\bu>\bo$,
\item[(ii)] if additionally $\dfrac{\b_{1}}{\L_{1}}+\dfrac{\b_{2}}{\L_{2}}<1$, there exists a
radial bound state $\bu^*>0.$ Furthermore,
$\Phi(\bu^*)>\Phi(\bv_{12}^*)=\Phi (\bv_1^*)+\Phi(\bv_2^*).$
\end{itemize}
\end{Theorem}
\begin{pf}
The proof of $(i)$ follows in a similar way as the one of Theorem
\ref{th:ground2} with appropriate changes, then we omit details for
short. With respect to $(ii)$, using Proposition \ref{prop:2}, one
can modify appropriately the arguments of the proof of Theorem
\ref{th:bound2} to obtain the result, so once again we can omit the
complete details.
\end{pf}
\begin{Remark}
It is easy to extend these results to systems with  any number of equations $N>3$ as the following
\be\label{eq:system-N}
\left\{\begin{array}{rcl}
-\D u +\l_0 u & = & \dyle u^3+\sum_{k=1}^{N-1}\b_{k}\, uv_k \\ & & \\
-\D v_j +\l_j v_j & = & \dyle\frac 12 v_j^2+\frac 12\b_{j} u^2; \qquad j=1,\cdots,N-1.\\
\end{array}\right.
\ee
 For example,  the existence of a positive radial ground state of  \eqref{eq:system-N} provided
$$
\b_{k}>\L_k=\inf_{\varphi\in H\setminus\{
0\}}\frac{\|\varphi\|_0^2}{\intR V_{k}^*\varphi^2dx}\qquad
k=1,\cdots N-1,
$$
where $V_k^*$ is the unique positive radial solution of $\D v+\l_k v=\frac 12 v^2$ in $E$, $k=1,\cdots, N-1$.
\end{Remark}
Another natural extension of \eqref{NLS-KdV2} to more than two
equations  different from \eqref{NLS-KdV2-3} is the following system
of NLS-NLS-KdV equations, \be\label{eq:NLS2-KdV}
\left\{\begin{array}{rcl}
-\D u_1 +\l_1 u_1 & = & u_1^3+\b_{12} u_1u_2^2+\b_{13}u_1v \\
-\D u_2 +\l_2 u_2 & = &  u_2^3+\frac 12\beta_{12} u_1^2u_2+\b_{23} u_2v\\
-\D v+\l v & = & \frac 12v^2+\frac 12\beta_{13} u_1^2+\frac
12\b_{23}u_2^2.
\end{array}\right.
\ee
Here we obtain a bifurcation result for this system in a similar
way as in Theorems \ref{th:2}, \ref{th:4}, that we enunciate in the
following.
\begin{Theorem}\label{th:44}
There exists $\e_0>0$ such that for any $0<\e <\e_0$ and $\b_{jk}=\e\wt{\b}_{jk}$, $k=1,2$, $j=2,3$, $k\neq j$,
System \eqref{eq:NLS2-KdV} has a  radial bound state $\bu_{\e}^*$ with  $\bu_{\e}\to \bu_0^*=(U_1,U_2,V)$ as $\e\to 0$. Moreover, if all $\b_{kj}>0$ then $\bu_{\e}^*>\bo$.
\end{Theorem}
 Where $U_j$  is the unique positive radial solution of $-\D u +\l_j u  = u^3$ in $E$; $j=1,\, 2$; and $V$ the corresponding positive radial solution to
 $-\D v +\l v = \frac 12 v^2$ in $E$.

Note that the non-negative radial semi-trivial solution $(0,0,V)$ is
a strict local minimum of the associated energy functional
constrained on the corresponding Nehari manifold provided
$$
\b_{j3}<\L_{j}=\inf_{\varphi\in H\setminus\{
0\}}\frac{\|\varphi\|_{\l_j}^2}{\intn V\varphi^2dx}\qquad j=1,2.
$$
While if either $\b_{13}>\L_1$ or $\b_{23}>\L_{2}$ then
 $(0,0,V)$ is a saddle critical point of $\Phi $ on $\cN$.

There also exist semi-trivial solutions coming from the solutions
studied in Section \ref{sec:mainR}, with the first or the
 second component $\equiv 0$. This fact makes different the analysis of
 \eqref{eq:NLS2-KdV} with respect to the previous studied systems \eqref{NLS-KdV2-3} and \eqref{eq:system-N}.

Finally, one could study more general extended systems of
\eqref{NLS-KdV2-3}, \eqref{eq:NLS2-KdV}
 with $N=m+\ell$; $m$-NLS and $\ell$-KdV coupled equations with $m,\,\ell\ge 2$ in the one dimensional case, or
 $N$-NLS equations if $n=2,3$. A careful analysis of this kind of systems including \eqref{eq:NLS2-KdV}
 will  be done in a forthcoming paper.
%%%%%%%%%%%%%%%%%%%%%%%%%%%%%%%%%


\begin{thebibliography}{99}
\bibitem{Ack} N. Akhmediev and A. Ankiewicz,  Solitons, Nonlinear pulses and beams, Champman \& Hall, London, 1997.

\bibitem{aa} J. Albert, J. Angulo Pava, {\it Existence and stability of ground-state solutions of a Schr\"odinger-KdV system}. Proc. Roy. Soc. Edinburgh Sect. A {\bf 133} (2003)
987-1029.

\bibitem{ab} J. Albert, S. Bhattarai, {\it Existence and stability of a two-parameter family of solitary waves for an NLS-KdV system}. Preprint, arXiv:1212.3902.

\bibitem{a} A. Ambrosetti, {\it A note on nonlinear Schr\"odinger systems: existence of a-symmetric solutions.} Adv. Nonlinear Stud. {\bf 6} (2006), no. 2, 149-155.

\bibitem{ac1} A. Ambrosetti, E. Colorado, {\it Bound and ground states of coupled nonlinear Schr\"odinger equations.} C. R. Math.
Acad. Sci. Paris {\bf 342} (2006), no. 7, 453-458.

\bibitem{ac2}  A. Ambrosetti, E. Colorado, {\it Standing waves of some coupled nonlinear Schr\"odin\-ger equations.} J. Lond. Math.
Soc. ({\bf 2}) 75 (2007), no. 1, 67-82.

\bibitem{acr} A. Ambrosetti, E. Colorado, D. Ruiz, {\it Multi-bump solitons to linearly coupled systems of nonlinear Schr\"odinger
equations.} Calc. Var. Partial Differential Equations {\bf 30} (2007), no. 1, 85-112.

\bibitem{a-m} A. Ambrosetti and A. Malchiodi, ``Perturbation methods and semilinear elliptic problems on $\Rn$", Progress in Math. Vol. {\bf 240}, Birkh\"auser, 2005.

\bibitem{ar}  A. Ambrosetti and P. H. Rabinowitz,   {\it Dual variational methods in critical point theory and applications.} J. Funct. Anal., {\bf 14} (1973), 349-381.

\bibitem{bt} T. Bartsch, Z.-Q. Wang, {\it Note on ground states of nonlinear Schr\"odinger systems.} J. Partial Differential
Equations {\bf 19} (2006), no. 3, 200-207.

\bibitem{chen-zou} Z. Chen; W. Zou, {\it An optimal constant for the existence of least energy solutions of a coupled Schr\"odinger system.} Calc. Var. Partial Differential Equations {\bf 48} (2013), no. 3-4, 695–711.

\bibitem{c-fract}  E. Colorado,  {\it Existence results for some systems of coupled fractional nonlinear Schr\"odinger equations.} Recent trends in nonlinear partial differential equations. II. Stationary problems, 135-150, Contemp. Math., {\bf 595}, Amer. Math. Soc., Providence, RI, 2013.

\bibitem{c}  E. Colorado, {\it Positive solutions to some systems of coupled nonlinear Schr\"odinger equations}.  Nonlinear Anal. {\bf 110} (2014) 104-112.

\bibitem{c2} E. Colorado, {\it Existence of Bound and Ground States for a System of Coupled Nonlinear Schr\"odinger-KdV
Equations}, Preprint, arXiv:1410.7638.

%\bibitem{ec-NLS-KdV} E. Colorado, {\it Positive solutions to a system of coupled nonlinear Schr\"odinger-KdV equations with general power nonlinearities}. In preparation.

\bibitem{cof} C.V. Coffman, {\it Uniqueness of the ground state solution for $\D u-u+u^3 = 0$ and a variational characterization of other solutions}, Arch. Rat. Mech. Anal. {\bf 46} (1972), 81-95.

\bibitem{cl} A.J. Corcho, F. Linares, {\it Well-posedness for the Schr\"odinger-Korteweg-de Vries system}. Trans. Amer. Math. Soc. {\bf 359} (2007) 4089-4106.

\bibitem{dfo} J.-P. Dias, M. Figueira, F. Oliveira, {\it Existence of bound states for the coupled Schr\"odinger-KdV system with cubic nonlinearity.} C. R. Math. Acad. Sci. Paris {\bf 348} (2010), no. 19-20, 1079-1082.

\bibitem{dfo2} J.-P. Dias, M. Figueira, F. Oliveira, {\it Well-posedness and existence of bound states for a coupled Schr\"odinger-gKdV system.} Nonlinear Anal. {\bf 73} (2010), no. 8, 2686-2698.

\bibitem{eke}  I. Ekeland, {\it On the variational principle}. J. Math. Anal.Appl. {\bf 47} (1974), 324-353.

\bibitem{fl} D. G. de Figueiredo, O. Lopes, {\it Solitary waves for some nonlinear Schrödinger systems.} Ann. Inst. H. Poincar\'e Anal. Non Lin\'eaire {\bf 25} (2008), no. 1, 149-161.

\bibitem{fo} M. Funakoshi, M. Oikawa, {\it The resonant Interaction between a Long Internal Gravity Wave and a Surface Gravity Wave Packet}. J. Phys. Soc. Japan. {\bf 52} (1983), no.1, 1982-1995.

\bibitem{iko} N. Ikoma, {\it Uniqueness of positive solutions for a nonlinear elliptic system.} NoDEA Nonlinear Differential Equations Appl. {\bf 16} (2009), no. 5, 555-567.

\bibitem{itanaka} N. Ikoma, K.  Tanaka, {\it A local mountain pass type result for a system of nonlinear Schr\"odinger equations.}
Calc. Var. Partial Differential Equations {\bf 40} (2011), no. 3-4, 449-480.

\bibitem{k} B. Kawohl, {\it Rearrangements and convexity of level sets in PDE}. Lecture Notes in Mathematics 1150 (Springer-Verlag, Berlin, 1985).

\bibitem{kw} M.K. Kwong, {\it Uniqueness of positive solutions of $\Delta u -u + u^p = 0$ in $\R^N$}. Arch. Rat. Mech. Anal.
\textbf{105} (1989), 243-266.

\bibitem{linwei} T.-C. Lin, J. Wei,  {\it Ground state of $N$ coupled nonlinear Schr\"odinger  equations in $\R^n$, $n\leq 3$}.
Comm. Math. Phys. {\bf 255} (2005), 629-653.

\bibitem{linwei2} T.-C. Lin, J. Wei, {\it Erratum: "Ground state of $N$ coupled nonlinear Schrödinger equations in ${\bf R}^n$, $n\leq 3$''} [Comm. Math. Phys. {\bf 255} (2005), no. 3, 629-653; MR2135447]. Comm. Math. Phys. {\bf 277} (2008), no. 2, 573-576.

\bibitem{Lions-JFA82} P.L. Lions,   {\em  Sym\'etrie et compacit\'e dans les espaces de Sobolev}. J. Funct. Anal., {\bf 49} (1982), no. 3, 315-334.

\bibitem{lions2} P.L. Lions, {\it The concentration-compactness principle in the calculus of variations. The locally compact case}. Ann.
Inst. H. Poincar\'e Anal. Non Lin\'eaire {\bf 1} (1984) 223-283.

\bibitem{liu-wang} Z. Liu, Z.-Q. Wang, {\it Ground states and bound states of a nonlinear Schr\"odinger system.} Adv. Nonlinear Stud. {\bf 10} (2010), no. 1, 175-193.

\bibitem{lz} C. Liu, Y. Zheng, {\it  On soliton solutions to a class of Schr\"odinger-KdV systems.} Proc. Amer. Math. Soc. {\bf 141} (2013), no. 10, 3477-3484.

\bibitem{mmp} L. Maia, E. Montefusco, B. Pellacci, {\it  Positive solutions for a weakly coupled nonlinear Schr\"odinger system.}
J. Differential Equations {\bf 229} (2006), no. 2, 743-767.

\bibitem{sirakov} B. Sirakov, {\it  Least energy solitary waves for a system of nonlinear Schr\"odinger equations in
$\Bbb R^n$.} Comm. Math. Phys. {\bf 271} (2007), no. 1, 199-221.

\bibitem{ww} J. Wei, T. Weth {\it Nonradial symmetric bound states for a system of coupled Schr\"odinger equations.}
Atti Accad. Naz. Lincei Cl. Sci. Fis. Mat. Natur. Rend. Lincei (9) Mat. Appl. {\bf 18} (2007), no. 3, 279-293.

\bibitem{wy} J. Wei, W. Yao, {\it Uniqueness of positive solutions to some coupled nonlinear Schr\"odinger equations.} Commun. Pure Appl. Anal. {\bf 11} (2012), no. 3, 1003-1011.
\end{thebibliography}
\end{document}